\documentclass{amsproc}

\newtheorem{theorem}{Theorem}[section]
\newtheorem{lemma}[theorem]{Lemma}
\newtheorem{proposition}[theorem]{Proposition}
\newtheorem{corollary}[theorem]{Corollary}

\theoremstyle{definition}
\newtheorem{definition}[theorem]{Definition}

\theoremstyle{remark}

\numberwithin{equation}{section}

\begin{document}

\title{Rank-one isometries of proper ${\rm CAT}(0)$-spaces}

\author{Ursula Hamenst\"adt}
\address{Mathematisches Institut der Universit\"at Bonn, 
Endenicher Allee 60, 53115 Bonn, Germany}
\email{ursula@math.uni-bonn.de}





\begin{abstract}
Let $X$ be a proper ${\rm CAT}(0)$-space with
visual boundary $\partial X$. Let $G$ be a
non-elementary group of isometries of $X$ with limit
set $\Lambda\subset \partial X$. 
We survey properties of the action of 
$G$ on $\Lambda$ under the assumption that 
$G$ contains a rank-one element. Among 
others, we show
that there is a dense orbit for the action of 
$G$ on the complement of the diagonal $\Delta$ in $\Lambda\times \Lambda$
and that pairs of 
fixed points of rank-one elements are dense in 
$\Lambda\times\Lambda-\Delta$.
\end{abstract}

\maketitle

\section{Introduction}

A geodesic metric space $(X,d)$ is called \emph{proper} if
closed balls in $X$ of finite radius are compact. 
A proper ${\rm CAT}(0)$-metric space $X$
can be compactified by adding the
\emph{visual boundary} $\partial X$.

The isometry group ${\rm Iso}(X)$ of a proper ${\rm CAT}(0)$-space 
$X$, equipped with the compact
open topology, is a locally compact $\sigma$-compact
topological group which acts as a group of homeomorphisms
on $\partial X$. The \emph{limit set} $\Lambda$ of a subgroup
$G$ of ${\rm Iso}(X)$ is the set of accumulation points
in $\partial X$ of an orbit of the action of $G$ on $X$. 
The limit set does not depend on the orbit, and 
it is closed and $G$-invariant.
The group $G$ is called
\emph{elementary} if either its limit set consists of at most
two points or if 
$G$ fixes a point in $\partial X$.

For every $g\in {\rm Iso}(X)$ the \emph{displacement function}
of $g$ is the function $x\to d(x,gx)$. The isometry $g$ is called
\emph{semisimple} if the displacement function assumes a
minimum on $X$. If this minimum vanishes
then $g$ has a fixed point in $X$ and is called \emph{elliptic},
and otherwise $g$ is called \emph{axial}.
If $g$ is axial then the subset $A$ of $X$ on which the displacement
function is minimal is isometric to a product 
space $C\times\mathbb{R}$ where
$C$ is a closed convex subset of $A$ and where $g$ acts on 
each of the geodesics $\{x\}\times \mathbb{R}$ as a translation.
Such a geodesic is called an \emph{axis} for $g$. We refer
to the books \cite{B95,BGS85,BH99} for basic properties of 
${\rm CAT}(0)$-spaces and for references.

Call an axial isometry $g$ of $X$ \emph{rank-one} if there
is an axis $\gamma$ for $g$ which does not bound a flat half-plane.
Here by a flat half-plane we mean 
a totally geodesic embedded isometric copy of an euclidean half-plane
in $X$.

Ballmann and Brin \cite{BB95} investigated 
discrete groups of isometries of a proper ${\rm CAT}(0)$-space
$X$ whose boundary $\partial X$ contains more than two points and
which act cocompactly on $X$. Such a group $G<{\rm Iso}(X)$
is necessarily non-elementary, and its limit set is the
whole boundary $\partial X$.
They showed (Theorem A of \cite{BB95})
that if $G$ 
contains a rank-one element then for any two 
non-empty open subsets $U,V$ of $\partial X$ there is 
an element $g\in G$ with $g(\partial X-U)\subset V$ and
$g^{-1}(\partial X-V)\subset U$. It is possible to choose
$g$ to be rank-one. Moreover (Theorem 4.6 of \cite{BB95}),
$G$ contains a free non-abelian subgroup.

In this note 
we extend this result 
to all non-elementary groups of isometries of a proper 
${\rm CAT}(0)$-space which contain a rank-one element.
For the formulation of our result, call the action of
a group $G$ of homeomorphism of a compact space \emph{minimal}
if every orbit is dense.

\begin{theorem}
Let $X$ be a proper ${\rm Cat}(0)$-space and let $G<{\rm Iso}(X)$
be a non-elementary subgroup with limit set $\Lambda$
which contains a rank-one element.
\begin{enumerate}
\item $\Lambda$ is perfect, and the $G$-action on $\Lambda$ is minimal.
\item Pairs of fixed points of rank-one elements are dense in 
the complement of the diagonal $\Delta$
of $\Lambda\times \Lambda$.
\item There is a dense orbit for the action of $G$ on
$\Lambda\times \Lambda-\Delta$.
\item $G$ contains a free subgroup with two generators
consisting of rank-one elements. 
\end{enumerate}
\end{theorem}

The second part of the above theorem is contained in 
\cite{CF08}. 
The last part of the above result is (in a slightly different context)
contained in \cite{BF08}, with a different proof.  
The paper \cite{CF08} also contains a proof of this last
part (without the conclusion that the free subgroup consists of
rank-one elements).
The other parts are probably also known to the experts, 
however they seem to be unavailable in the literature in this form.

The organization of this paper is as follows. In Section 2 we 
collect some basic geometric properties of proper ${\rm CAT}(0)$-spaces.
In Section 3 we look at 
geodesics in proper ${\rm CAT}(0)$-spaces.
Following \cite{BF08}, we define contracting geodesics and study
some of their properties. 
In Section 4 we investigate some
geometric properties of isometries of a
proper ${\rm CAT}(0)$-space $X$, and in Section 5 we
look at groups of isometries and prove the theorem above.

\section{Basic ${\rm CAT}(0)$-geometry}

In this section we summarize some geometric properties
of ${\rm CAT}(0)$-spaces. We use
the books \cite{B95,BGS85,BH99} as our main references and for
the discussion of a large set of examples.

A \emph{triangle} $\Delta$ in a geodesic metric space consists
of three vertices connected by three (minimal) geodesic arcs $a,b,c$. 
A comparison triangle $\bar \Delta$ for $\Delta$ 
in the euclidean plane
is a triangle in $\mathbb{R}^2$ 
with the same side-lengths as $\Delta$. By the triangle inequality,
such a comparison triangle exists always, and it is 
unique up to isometry. For a point $x\in a\subset \Delta$
the comparison point of $x$ 
in the comparison triangle $\bar\Delta$ 
is the point on the side $\bar a$ of $\bar\Delta$ corresponding to $a$
whose distance to the endpoints of $\bar a $
coincides with the distance of $x$ to the
corresponding endpoints of $a$.

\begin{definition}\label{cat}
A geodesic metric space $(X,d)$ 
is called a \emph{${\rm CAT}(0)$-space}
if for every geodesic triangle $\Delta$ in $X$ with 
sides $a,b,c$ and every comparison triangle
$\bar\Delta$ 
in the euclidean plane 
with sides $\bar a,\bar b,\bar c$  
and for all $x,y\in \Delta$ and all 
comparison points 
$\bar x,\bar y\in \bar \Delta$ we have
\[d(x,y)\leq d(\bar x,\bar y).\]
\end{definition}

A complete ${\rm CAT}(0)$-space is called a 
\emph{Hadamard space}. 
In a Hadamard space $X$, any two points can
be connected by a unique geodesic which varies
continuously with the endpoints. The distance
function is convex: If $\gamma,\zeta$ are two
geodesics in $X$ parametrized on the same interval 
then the function 
$t\to d(\gamma(t),\zeta(t))$ is convex.
More generally, we call a function
$f:X\to \mathbb{R}$ \emph{convex} if
for every geodesic $\gamma:J\to \mathbb{R}$ the
function $t\to f(\gamma(t))$ is convex \cite{B95}.

For a fixed point $x\in X$,
the \emph{visual boundary} $\partial X$ of $X$
is defined to be the space
of all geodesic rays issuing from $x$ equipped with the
topology of uniform convergence on bounded sets.
This definition is independent of the choice of $x$.
We denote the point in $\partial X$ defined by 
a geodesic ray $\gamma:[0,\infty)\to X$ by $\gamma(\infty)$.
We also say that $\gamma$ 
\emph{connects} $x$ to 
$\gamma(\infty)$. If $X$ is proper then the
visual boundary of $X$ is compact \cite{B95}.
Note that if $\gamma_1,\gamma_2:[0,\infty)\to X$
are two geodesic rays with 
$\gamma_1(\infty)=\gamma_2(\infty)$ then the
function $t\to d(\gamma_1(t),\gamma_2(t))$ is convex and 
bounded (see Chapter II.8 of \cite{BH99}) 
and hence it is nonincreasing.

There is another description of the visual boundary of $X$ 
as follows.
Let $C(X)$ be the space of all continuous functions on $X$
endowed with the topology of uniform convergence on
bounded sets. 
Fix a point $y\in X$ and for $x,z\in X$ define
\[b_x(y,z)=d(x,z)-d(x,y).\]
Then we have 
\begin{equation}\label{sym}
b_x(y,z)=-b_x(z,y) \text{ for all } y,z\in X
\end{equation}
and 
\begin{equation}\label{lip}
\vert b_x(y,z)-b_x(y,z^\prime)\vert \leq d(z,z^\prime)
\text{ for all }z,z^\prime\in X\end{equation}
and hence the function 
$b_{x}(y,\cdot):z\to b_x(y,z)$ is 
one-Lipschitz and vanishes at $y$.
The assignment $x\to b_x(y,\cdot)$ is an embedding of $X$ into
$C(X)$. Moreover, for every $x\in X$ the function
$b_x(y,\cdot)$ is convex. If $\tilde y\in X$ is 
another basepoint then we have 
\begin{equation}\label{add}
b_x(\tilde y,\cdot)=b_x(y,\cdot)+b_x(\tilde y,y).\end{equation}

A sequence $\{x_n\}\subset X$ \emph{converges at infinity}
if $d(x_n,x)\to \infty$ and if the functions 
$b_{x_n}(y,\cdot)$ converge in $C(X)$. The visual boundary 
$\partial X$ of $X$ can also be defined as the 
subset of $C(X)$ of all functions 
which are obtained
as limits of functions $b_{x_n}(y,\cdot)$ for
sequences $\{x_n\}\subset X$ 
which converge at infinity. In particular, 
the union $X\cup \partial X$
is naturally a closed subset of $C(X)$.

Namely,
if $\gamma:[0,\infty)\to X$ is any geodesic ray then 
for every sequence $t_n\to \infty$ the sequence $\{\gamma(t_n)\}$
converges at infinity, and the limit function 
$b_{\gamma(\infty)}(y,\cdot)$ does not depend on 
the sequence $\{t_n\}$. The 
function $b_{\gamma(\infty)}(y,\cdot)$ 
is called a \emph{Busemann function} at 
$\gamma(\infty)\in \partial X$.
The Busemann function
$b_{\gamma(\infty)}(\gamma(0),\cdot)$ satisfies
\begin{equation}\label{buseray}
b_{\gamma(\infty)}(\gamma(0),\gamma(t))=-t\text{ for all }t\geq 0,
\end{equation}
moreover it is convex.

Vice versa, if the sequence
$\{x_n\}\subset X$ converges at infinity then 
the geodesics $\gamma_n$ connecting a fixed point $x\in X$ to $x_n$
converge locally uniformly to a geodesic ray which only
depends on the limit of the functions $b_{x_n}(y,\cdot)$
(Chapter II.1 and II.2 of \cite{B95}).

From now on let $X$ be a \emph{proper} (i.e. complete and 
locally compact) ${\rm CAT}(0)$-space.
Then $X\cup \partial X$ is compact. A subset $C\subset X$ is 
\emph{convex} if for $x,y\in C$ the geodesic connecting
$x$ to $y$ is contained in $C$ as well. 
For every closed convex set $C\subset X$
and every $x\in X$ there is a unique point
$\pi_C(x)\in C$ of smallest distance to $x$ (Proposition II.2.4 of 
\cite{BH99}).
Now let $J\subset \mathbb{R}$ be a closed connected set and
let $\gamma:J\to X$ be a geodesic arc. Then $\gamma(J)\subset X$
is closed and convex and hence there is a 
shortest distance projection $\pi_{\gamma(J)}:X\to \gamma(J)$.
Then $\pi_{\gamma(J)}(x)$ is the unique minimum for the
restriction of the function $b_x(y,\cdot)$ to $\gamma(J)$.
By equality (\ref{add}), 
this does not depend on the choice of the basepoint $y\in X$.
The projection $\pi_{\gamma(J)}:X\to \gamma(J)$ is 
distance non-increasing.

For $\xi\in \partial X$
the function $t\to b_\xi(y,\gamma(t))$
is convex. Let $\overline{\gamma(J)}$ be the closure of $\gamma(J)$
in $X\cup \partial X$. If $b_\xi(y,\cdot)\vert \gamma(J)$ 
assumes a
minimum then we can define $\pi_{\gamma(J)}(\xi)\subset \overline{\gamma(J)}$ 
to be the closure in $\overline{\gamma(J)}$ of the  
connected subset of $\gamma(J)$ consisting of all such minima. 
If $b_\xi(y,\cdot)\vert \gamma(J)$ does not assume a minimum
then by continuity the set $J$ is unbounded and by convexity 
either $\lim_{t\to \infty}b_\xi(y,\gamma(t))=
\inf\{b_\xi(y,\gamma(s))\mid s\in J\}$ 
or $\lim_{t\to -\infty}b_\xi(y,\gamma(t))=
\inf\{b_\xi(y,\gamma(s))\mid s\in J\}$. In the first
case we define 
$\pi_{\gamma(J)}(\xi)=\gamma(\infty)\in \partial X$, and in the
second case we define  
$\pi_{\gamma(J)}(\xi)=\gamma(-\infty)$. Then for every $\xi\in \partial X$
the set  
$\pi_{\gamma(J)}(\xi)$ is a closed connected 
subset of $\overline{\gamma(J)}$ (which may
contain points in both $X$ and $\partial X$).

The following simple observation will be useful several
times in the sequel.

\begin{lemma}\label{projectionfixed}
Let $\gamma:J\to X$ be a geodesic,
let $\xi\in \partial X$ and assume that 
$\pi_{\gamma(J)}(\xi)\cap X\not=\emptyset$. 
If $c:[0,\infty)\to X$ is a geodesic ray connecting
a point $c(0)\in \pi_{\gamma(J)}(\xi)$ 
to $c(\infty)=\xi$ then
$\pi_{\gamma(J)}(c(t))=c(0)$ for all $t\geq 0$.
\end{lemma}
\begin{proof}
If $c:[0,\infty)\to X$ connects 
$c(0)\in \pi_{\gamma(J)}(\xi)$
to $\xi$ then property (\ref{buseray}) of 
Busemann functions implies that $b_\xi(c(0),c(R))
=-b_\xi(c(R),c(0))=-R$ for all $R>0$.
Moreover, $b_\xi(c(R),\cdot)$ is one-Lipschitz and hence
if $z=\pi_{\gamma(J)}(c(R))\not=c(0)$ then
$d(z,c(R))<R$ and consequently 
$b_\xi(c(R),z)<R$. However, this implies that 
\[b_\xi(c(0),z)=b_\xi(c(0),c(R))+b_\xi(c(R),z)<0\] 
which violates the
assumption that $c(0)\in \pi_{\gamma(J)}(\xi)$.
The lemma is proven.
\end{proof}

We also note the following easy fact.

\begin{lemma}\label{contrcont}
Let $\gamma:J\to X$ be a geodesic and let 
$(x_i)\subset X$ be a sequence converging to some
$\xi\in \partial X$. Then up to passing 
to a subsequence, the sequence $\pi_{\gamma(J)}(x_i)$ 
converges to a point in $\pi_{\gamma(J)}(\xi)$.
\end{lemma}
\begin{proof}
Since the closure $\overline{\gamma(J)}$ of $\gamma(J)$
in $X\cup \partial X$ is compact, up to passing to 
a subsequence the sequence $\pi_{\gamma(J)}(x_i)$ converges to 
a point $z\in \overline{\gamma(J)}$. On the other hand,
the functions $b_{x_i}(x,\cdot)$ converge as
$i\to \infty$ locally uniformly to the Busemann function
$b_\xi$. Now if $z\in \gamma(J)$ then this implies that
$z$ is a minimum for $b_\xi(x,\cdot)$ and hence
$z\in \pi_{\gamma(J)}(\xi)$ by definition. 

Otherwise assume that $J\supset [a,\infty)$ for some $a\in 
\mathbb{R}$ and that $\pi_{\gamma(J)}(x_i)=\gamma(t_i)$ 
where $t_i\to \infty$. Then by convexity,  
for every $s\in (a,\infty)$, every $t>s$ 
and every $i$ which is sufficiently large that $t_i>t$ we have  
$b_{x_i}(x,\gamma(t))\leq b_{x_i}(x,\gamma(s))$. 
Since $b_{x_i}(x,\cdot)\to b_\xi(x,\cdot)$ locally uniformly, 
we also have $b_\xi(x,\gamma(t))\leq b_\xi(x,\gamma(s))$ and hence
indeed $\lim_{t\to \infty}
b_\xi(x,\gamma(t))=\inf\{b_\xi(x,\gamma(s))\mid s\in J\}.$
This shows the lemma.
\end{proof}

\section{Contracting geodesics}

In this section we discuss some geometric properties 
of geodesics in a proper ${\rm CAT}(0)$-space $X$. As a convention,
geodesics are always defined on closed connected
subsets of $\mathbb{R}$.
We begin with the following definition which is due to Bestvina and
Fujiwara (Definition 3.1 of \cite{BF08}).

\begin{definition}\label{contracting}
A geodesic arc $\gamma:J\to X$ is 
\emph{$B$-contracting}
for some $B>0$ if for every closed 
metric ball $K$ in $X$ which is disjoint
from $\gamma(J)$ the diameter of the 
projection $\pi_{\gamma(J)}(K)$ does not exceed $B$. 
\end{definition}
We call a geodesic \emph{contracting} if it is $B$-contracting
for some $B>0$.
As an example, every geodesic in a ${\rm CAT}(\kappa)$-space for some
$\kappa<0$ is $B$-contracting for a number $B=B(\kappa)>0$ only
depending on $\kappa$.

The next lemma (Lemma 3.2 and 3.5 of \cite{BF08}) shows that 
a connected subarc of a contracting geodesic is contracting and that
a triangle containing a $B$-contracting
geodesic as one of its sides is
uniformly thin.

\begin{lemma}\label{thintriangle} 
\begin{enumerate}
\item Let $\gamma:J\to X$ be a $B$-contracting
geodesic. Then for every closed connected subset $I\subset J$,
the subarc $\gamma(I)$ of $\gamma$ is $B+3$-contracting.
\item
Let $\gamma:[a,b]\to X$ be a $B$-contracting geodesic.
If $x\in X$ is such that $\pi_{\gamma[a,b]}(x)=a$ then for every
$t\in [a,b]$ the geodesic connecting $x$ to $\gamma(t)$
passes through the $3B+1$-neighborhood of $\gamma(a)$. 
\end{enumerate}
\end{lemma}

Lemma \ref{thintriangle} implies that for 
a $B$-contracting biinfinite geodesic $\gamma:\mathbb{R}\to 
X$ and for 
$\xi\in \partial X-\{\gamma(-\infty),
\gamma(\infty)\}$ the projection $\pi_{\gamma(\mathbb{R})}(\xi)$ 
is a compact subset of $\gamma(\mathbb{R})$ of diameter at most
$6B+4$.

\begin{lemma}\label{projectiondiameter}
For some $B>0$ 
let $\gamma:\mathbb{R}\to X$ be a biinfinite 
$B$-contracting geodesic. Then for $\xi\in 
\partial X-\{\gamma(-\infty),\gamma(\infty)\}$ 
the restriction to $\gamma(\mathbb{R})$ of a Busemann function 
$b_\xi$ at $\xi$ is bounded from below and
assumes a minimum. If $T\in \mathbb{R}$ is
such that $\gamma(T)\in 
\pi_{\gamma(\mathbb{R})}(\xi)$ then 
\[\vert T-t\vert-3B-2 \leq b_\xi(\gamma(T),\gamma(t))\leq \vert T-t\vert
\text{ for all }t\in \mathbb{R}.\] 
\end{lemma}
\begin{proof} Let $\gamma:\mathbb{R}\to X$ be
a biinfinite $B$-contracting geodesic
and let $\xi\in \partial X-\{\gamma(\infty),\gamma(-\infty)\}$.
Let $c:[0,\infty)\to X$ be the geodesic 
ray which connects $\gamma(0)$ to $\xi$. For $R>0$ let
$t_R\in \mathbb{R}$ be such that $\pi_{\gamma(\mathbb{R})}(c(R))=
\gamma(t_R)$. By Lemma \ref{thintriangle}, the
geodesic $c$ passes through the $3B+1$-neighborhood of 
$\gamma(t_R)$. By triangle comparison, the geodesic segment
$\gamma[0,t_R]$ is contained in the $3B+1$-neighborhood of
$c[0,\infty)$. Thus if there is a sequence
$R_i\to \infty$ such that $t_{R_i}\to \infty$
(or $t_{R_i}\to -\infty)$ then 
the geodesic ray $\gamma[0,\infty)$ 
(or $\gamma(-\infty,0])$ 
is contained in the $3B+1$-neighborhood
of $c[0,\infty)$ and hence $c=\gamma[0,\infty)$ 
(or $c=\gamma(-\infty,0]$) which
is impossible.  

As a consequence, the set $\{t_R\mid R\geq 0\}\subset
\mathbb{R}$ is bounded
and therefore there is a number $T\in \mathbb{R}$ and a 
sequence $R_i\to \infty$
such that $t_{R_i}\to T$ $(i\to \infty)$. 
By Lemma \ref{thintriangle}, for sufficiently large $i$ 
and all $t\in \mathbb{R}$ the geodesic connecting $\gamma(t)$ to 
$c(R_i)$ passes through the $3B+2$-neighborhood of
$\gamma(T)$. On the other hand, as $i\to \infty$ these
geodesics converge locally uniformly to the geodesic connecting 
$\gamma(t)$ to $\xi$. 
Together with (\ref{sym}), (\ref{lip}) and (\ref{buseray}) above,
this implies that 
\[\vert T-t\vert -3B-2 
\leq b_\xi(\gamma(T),\gamma(t))\leq \vert T-t\vert\]
as claimed in the lemma. In particular,     
the restriction of the function $b_\xi(\gamma(T),\cdot)$ to 
$\gamma(\mathbb{R})$ is
bounded from below by $-3B-2$, and if $\vert T-t\vert >3B+2$ then 
$b_\xi(\gamma(T),\gamma(t))>0$ and hence
$\gamma(t)\not\in \pi_{\gamma(\mathbb{R})}(\xi)$.
\end{proof}

{\bf Remark:} Lemma \ref{thintriangle} and Lemma \ref{projectiondiameter}
and their proofs are valid without the assumption that the space $X$
is proper.

\bigskip

Two points $\xi\not=\eta\in \partial X$ 
are connected in $X$ by a geodesic
if there is a geodesic $\gamma:\mathbb{R}\to X$ 
with $\gamma(\infty)=\xi,\gamma(-\infty)=\eta$. 
Unlike in a proper hyperbolic geodesic metric space,
such a geodesic need not exist. Therefore we define.

\begin{definition}\label{visibility}
A point $\xi\in \partial X$ is called a \emph{visibility
point} if for every $\zeta\not=\xi\in \partial X$ there
is a geodesic connecting $\xi$ to $\zeta$.
\end{definition}

Lemma \ref{thintriangle} and Lemma \ref{projectiondiameter}
are used to show (compare also Lemma 23 of \cite{K05}).

\begin{lemma}\label{visibilitypoint}
Let $\gamma:[0,\infty)\to X$ be a contracting geodesic ray. Then 
$\gamma(\infty)\in \partial X$ 
is a visibility point.
\end{lemma}
\begin{proof}
Let $\gamma:[0,\infty)\to X$ be a $B$-contracting geodesic ray 
for some $B>0$ and let
$\xi\in \partial X-\gamma(\infty)$.
By Lemma \ref{projectiondiameter} (or, rather, its obvious 
modification for geodesic rays)  
the projection
$\pi_{\gamma[0,\infty)}(\xi)$ is a compact subset of $\gamma[0,\infty)$ of 
diameter at most $6B+4$. 

Let $r\geq 0$ be such that
$\gamma(r)\in \pi_{\gamma[0,\infty)}(\xi)$.
Let $c:[0,\infty)\to X$ be the geodesic ray connecting
$c(0)=\gamma(r)$ to $\xi$. By Lemma \ref{projectionfixed}, 
for every $t>0$ we have
$\pi_{\gamma[0,\infty)}(c(t))=\gamma(r)$. 
Together with  
Lemma \ref{thintriangle}, this shows that 
for every $t>0$ the geodesic $\zeta_t$ 
connecting
$\gamma(t)$ to $c(t)$ passes through 
the $3B+1$-neighborhood of $\gamma(r)$. Since
$X$ is proper, up to reparametrization and up to
passing to a subsequence we may assume that the geodesics
$\zeta_t$ converge uniformly on compact sets as $t\to \infty$ to a
geodesic $\zeta$. By construction and by convexity, 
$\zeta$ connects $\gamma(\infty)$ to $c(\infty)=\xi$. 
Since $\xi\in \partial X-\gamma(\infty)$
was arbitrary, this shows the lemma.
\end{proof}

For fixed $B>0$, $B$-contracting geodesics are stable under
limits.

\begin{lemma}\label{bconvergence}
Let $B>0$ and let $\gamma_i:J_i\to X$ be a
sequence of $B$-contracting geodesics converging locally
uniformly to a geodesic $\gamma:J\to X$. Then 
$\gamma$ is $B$-contracting.
\end{lemma}
\begin{proof}
Assume to the contrary that
there is a sequence $(\gamma_i:J_i\to X)$ of $B$-contracting
geodesics in $X$ converging locally uniformly to a geodesic
$\gamma:J\to X$ which is not $B$-contracting. Then there is a compact
metric ball $K$ which is disjoint from $\gamma(J)$ and such that
the diameter of $\pi_{\gamma(J)}(K)$ is bigger than $B$.
In other words, there are two points $x,y\in K$ with
$d(\pi_{\gamma(J)}(x),\pi_{\gamma(J)}(y))>B$.

Since $\gamma_i\to \gamma$ locally uniformly, for sufficiently
large $i$ the ball $K$ is disjoint from $\gamma_i$. 
Let $u_i=\pi_{\gamma_i(J_i)}(x),z_i=\pi_{\gamma_i(J_i)}(y)$.
If $i>0$ is sufficiently large that $K$ is disjoint from $\gamma_i$ then
we have $d(u_i,z_i)\leq B$ since $\gamma_i$ is $B$-contracting. Moreover,  
the distance between $u_i$ and $x$ and between $z_i$ and $y$
is bounded independently of $i$. Thus
up to passing to a subsequence we may assume that
$u_i\to u,z_i\to z$. Then $u,z\in \gamma(J)$ and
$d(u,z)\leq B$ by continuity and therefore 
up to possibly exchanging $x$ and $y$ we may
assume that $u\not=\pi_{\gamma(J)}(x)$. Since the shortest
distance projection of $x$ into $\gamma(J)$ is unique, 
we have $d(u,x)>d(\pi_{\gamma(J)}(x),x)$.
But $\gamma_i\to \gamma$ locally uniformly and $d(u_i,x)\to d(u,x)$ and
therefore for sufficiently large $i$ the point 
$\pi_{\gamma_i(J_i)}(\pi_{\gamma(J)}(x))\in\gamma_i(J_i)$
is closer to $x$ than $u_i$. This contradicts the choice of $u_i$ and
shows the lemma.
\end{proof}

\section{Rank-one isometries}

As before, let $X$ be a proper ${\rm CAT}(0)$-space.
For an isometry $g$ of $X$ define the \emph{displacement function}
$d_g$ of $g$ to be the function $x\to d_g(x)=d(x,gx)$.

\begin{definition}\label{semisimple}
An isometry $g$ of $X$ 
is called \emph{semisimple} if $d_g$ achieves
its minimum in $X$. If $g$ is semisimple and ${\rm min}\,d_g=0$
then $g$ is called \emph{elliptic}. A semisimple isometry
$g$ with ${\rm min}\,d_g>0$ is called \emph{axial}.
\end{definition}

By the above definition, an isometry is elliptic if and only
if it fixes at least one point in $X$. Any isometry of $X$
which admits a bounded orbit in $X$ is elliptic \cite{BH99}.
By Proposition II.3.3 of \cite{B95}, an isometry $g$ of $X$ is
axial if and only if there is a geodesic
$\gamma:\mathbb{R}\to X$ such that
$g\gamma(t)=\gamma(t+\tau)$ for every $t\in \mathbb{R}$ where
$\tau=\min\,d_g>0$.
Such a geodesic is called an \emph{oriented axis} for $g$.
Note that the geodesic $t\to \gamma(-t)$ is an oriented 
axis for $g^{-1}$. The endpoint $\gamma(\infty)$ of $\gamma$
is a fixed point for the action of $g$ on $\partial X$
which is called the \emph{attracting fixed point}. 
The closed convex set $A\subset X$ of all points for which the displacement
function of $g$ is minimal is 
isometric to $C\times \mathbb{R}$ where
$C\subset A$ is closed and convex. For each
$x\in C$ the set $\{x\}\times \mathbb{R}$ is an axis of $g$.

The following definition is due to Bestvina
and Fujiwara (Definition 5.1 of \cite{BF08}).

\begin{definition}\label{rankone}
An isometry $g\in {\rm Iso}(X)$ is called 
\emph{$B$-rank-one} for some $B>0$ 
if $g$ is axial and admits a $B$-contracting axis.
\end{definition}
We call an isometry $g$ \emph{rank-one} if $g$ is $B$-rank-one for some
$B>0$.

The following statement is Theorem 5.4 of \cite{BF08}.

\begin{proposition}\label{flatrank}
An axial isometry of $X$ with axis $\gamma$ is rank-one
if and only if $\gamma$ does not bound a flat half-plane.
\end{proposition}

A homeomorphism $g$ of a compact space $K$ is said
to act with \emph{north-south dynamics} if there are
two fixed points $a\not=b\in K$ for the action of $g$ such that
for every neighborhood $U$ of $a$, $V$ of $b$ there is some 
$k>0$ such that $g^k(K-V)\subset U$ and 
$g^{-k}(K-U)\subset V$.
The point $a$ is called the \emph{attracting fixed point} for $g$,
and $b$ is the \emph{repelling fixed point}.
A rank-one isometry acts with north-south dynamics on $\partial X$
(see Lemma 3.3.3 of \cite{B95}). For completeness of exposition,
we provide a proof of this fact.

\begin{lemma}\label{northsouth}
An axial isometry $g$ of $X$ is rank-one if and only if
$g$ acts with north-south dynamics
on $\partial X$.
\end{lemma}
\begin{proof} Let $g$ be a $B$-rank-one isometry of $X$ and let 
$\gamma$ be a $B$-contracting axis of $g$.
Let
$A\subset \gamma$ be a compact connected fundamental domain 
for the action of $g$ on 
$\gamma$. A Busemann function
$b_\xi(x,\cdot)$ depends continuously on $\xi\in \partial X$,
and $b_{\gamma(\infty)}(\gamma(0),\gamma(t))=-t$ and
$b_{\gamma(-\infty)}(\gamma(0),\gamma(t))=t$ for all $t\in \mathbb{R}$.
By Lemma \ref{projectiondiameter}, 
for every $\xi\in \partial X-\{\gamma(\infty),\gamma(-\infty)\}$ the set
$\pi_{\gamma(\mathbb{R})}(\xi)$ is a compact subset of $\gamma(\mathbb{R})$
of diameter at most $6B+4$. Therefore the set  
$K=\{\xi\in  \partial X-\{\gamma(\infty),
\gamma(-\infty)\}\mid \pi_{\gamma(\mathbb{R})}(\xi)\cap A\not=\emptyset\}$
is closed and does not contain $\gamma(\infty),\gamma(-\infty)$.

Using again Lemma \ref{projectiondiameter}, 
for every $\xi\in \partial X-\{\gamma(\infty),\gamma(-\infty)\}$ 
there is some 
$k\in \mathbb{Z}$ 
with $\pi_{\gamma(\mathbb{R})}(\xi)\cap  g^kA\not=\emptyset$.
By equivariance under the action of the infinite cyclic subgroup of
${\rm Iso}(X)$ generated by $g$, this means that $\xi\in g^kK$ and hence
$\cup_kg^k(K)=\partial X-\{\gamma(-\infty),\gamma(\infty)\}$. 

By comparison, for every 
neighborhood $V$ of $\gamma(-\infty)$ there is a number 
$m>0$ such that $V$ contains the endpoints of all 
geodesic rays $\zeta:[0,\infty)\to \partial X$ 
issuing from $\gamma(0)=\zeta(0)$
which pass through the $3B+1$-neighborhood of $\gamma(-t)$
for some $t\geq m$.
Lemma \ref{thintriangle} and  
Lemma \ref{projectionfixed} show 
that if $\xi\in \partial X$ is 
such that there is a point
$z\in \pi_{\gamma(\mathbb{R})}(\xi)\cap \gamma(-\infty,-m]$
then the geodesic ray connecting $\gamma(0)$ to
$\xi$ intersects the $3B+1$-neighborhood of $z$
and hence $\xi\in V$.
This means that if $k\geq 0$ is sufficiently large that
$g^{-k}A\subset \gamma(-\infty,-m]$ then 
$\cup_{j\leq -k}g^{j}K
\subset V$ and hence 
$\partial X-V\subset
\cup_{j> -k}g^jK\cup\{\gamma(\infty)\}$.
Similarly, for every neighborhood $U$ of $a$ 
there is some $\ell >0$ such that $\cup_{j\geq \ell}g^j K\subset U$.
Then by equivariance, we have 
$g^{k+\ell}(\partial X-V)\subset U$.
This shows that $g$ acts with north-south dynamics
on $\partial X$.

On the other hand, let $g$ be an axial isometry
of $X$ which acts with north-south dynamics on $\partial X$, with
attracting and repelling fixed point $a,b\in \partial X$,
respectively.  
If $g$ is not rank-one then there is an oriented axis 
$\gamma$ for $g$ which
bounds a flat half-plane $F\subset X$. 
For $k\in \mathbb{Z}$ the 
image $g^kF$ of $F$ under the isometry
$g^k$ is again a flat
half-plane with boundary $\gamma$.

By the definition
of the topology on $\partial X$, there is a 
neighborhood $U$ of $a$
in $\partial X$ with the following property.
Let $\xi:[0,\infty)\to X$ be a geodesic ray
with $\xi(0)=\gamma(0)$ which
is contained in a flat half-plane $G$ bounded by $\gamma$.
If $\xi$ encloses an angle with
the oriented ray $\gamma[0,\infty)$ 
in $G$ which is bigger than $\pi/4$ 
then $\xi(\infty)\not\in U$. 

Let $\xi:[0,\infty)\to F$ be the geodesic ray with $\xi(0)=\gamma(0)$ 
which meets 
$\gamma(\mathbb{R})$ perpendicularly at $\gamma(0)$. Since by assumption
$g$ acts with north-south dynamics on $\partial X$, 
there is some $k>0$ such that $g^k\xi(\infty)\in U$. On the other
hand, $g^k\xi$ is a geodesic
ray in the flat half-plane $g^kF$ which
is perpendicular to $\gamma$.
Then the angle in $g^kF$ 
between $\gamma$ and the ray in $g^kF$ 
which connects $x$ to $g^k\xi(\infty)$ equals $\pi/2$ as well.
However, this means that $g^k\xi(\infty)\not\in U$ which is
a contradiction. Therefore $g$ is rank-one. The lemma is proven.
\end{proof} 

{\bf Remark:} The proof of the statement that
a rank-one isometry of $X$ acts with 
north-south dynamics on $\partial X$ is valid
without the assumption that $X$ is proper.

\bigskip

We conclude this section with a characterization of
rank-one isometries which is easier to verify.

\begin{lemma}\label{northsouth2}
Let $g\in {\rm Iso}(X)$ and assume 
that there are non-trivial open 
subsets $V_1,V_2$ of $\partial X$ with 
the following properties.
\begin{enumerate}
\item The closures 
$\overline{V_1},\overline{V_2}$ of $V_1,V_2$ are disjoint.
\item There is a number $B>0$ and there is a $B$-contracting 
biinfinite geodesic with both endpoints in 
$\partial X-\overline{V_1}-\overline{V_2}$.
\item The distance between any biinfinite $B$-contracting geodesic
with both
endpoints in $V_1$ 
and any biinfinite $B$-contracting 
geodesic with both endpoints
in $V_2$ 
is bounded from below by a universal positive constant.
\item
$g(\partial X-V_2)\subset V_1$ and 
$g^{-1}(\partial X-V_1)\subset V_2$. 
\end{enumerate}
Then 
$g$ is rank-one, with attracting fixed point in $V_1$ and 
repelling fixed point in $V_2$.
\end{lemma}
\begin{proof}
Let $g\in {\rm Iso}(X)$ 
be any isometry with the properties
stated in the lemma for
open subsets $V_1,V_2$ of $\partial X$ with disjoint
closure $\overline{V_1},\overline{V_2}$.
We have to show that $g$ is rank-one. Note
that $g$ maps a nontrivial open
neighborhood of $\overline{V_1}$ into
$\overline{V_1}$ and hence the order of $g$ is infinite. 
We show first that  
$g$ is not elliptic. 

For this assume to the contrary that $g$ is elliptic.
Then every orbit in $X$ of the infinite cyclic subgroup $G$ of
${\rm Iso}(X)$ generated by $g$ is bounded. 
Let $\gamma$ be a biinfinite $B$-contracting geodesic
whose endpoints $\gamma(\infty),\gamma(-\infty)$ are
contained in $\partial X-\overline{V_1}-\overline{V_2}$.
Such a geodesic exists by the
second requirement for $V_1,V_2$ stated in the lemma.
For every $k\in \mathbb{Z}$ the geodesic $g^k\gamma$ is $B$-contracting. 

Since the $G$-orbit of every point in $X$ is bounded,
the geodesics $g^k\gamma$ all pass through a fixed compact subset
of $X$. By the Arzela Ascoli theorem,
for every sequence $k_i\to \infty$ there is a subsequence 
$(k_{i_j})$ such that 
the geodesics $g^{k_{i_j}}\gamma$ converge locally uniformly
to a biinfinite geodesic in $X$. Lemma \ref{bconvergence}
shows that such a limiting geodesic 
$\zeta$ is $B$-contracting. Moreover, by the properties of $g$
the endpoints of $\zeta$ in $\partial X$ 
are contained in $\overline{V_1}$.
The collection 
${\mathcal G}_+$ of all such limiting geodesics
is closed for the compact open topology, moreover
${\mathcal G}_+$ is $g$-invariant. 

Let $x\in X$ be a fixed point of $g$ and let 
$K_+=\{\pi_{\zeta(\mathbb{R})}(x)\mid
\zeta\in {\mathcal G}_+\}$. Since 
${\mathcal G}_{+}$ is closed and $g$-invariant 
and consists of $B$-contracting geodesics passing through
a fixed compact subset of $X$, the set $K_+$ is
compact and $g$-invariant. Moreover, 
by the third requirement in the statement of the lemma,
there is a number $c>0$ such that 
the distance in $X$ of $K_+$ to any biinfinite $B$-contracting
geodesic with both endpoints in $\overline{V_2}$ is at least $2c$. 
On the other hand, by definition of the set $K_+$, 
for $y\in K_+$ there is some $z\in X$ with $d(y,z)<c$ and there is 
some $\ell >0$ such that $g^{-\ell}(z)\in \gamma$.
Then for $m>\ell$, the point $g^{-m}(z)$ lies on a $B$-contracting geodesic
connecting two points in $\overline{V_2}$. In particular, its distance
to $K_+$ is at least $2c$. 
Since $K_+$ is $g$-invariant, this implies that
$d(g^{-m}(y),g^{-m}(z))\geq 2c$, on the other hand also
$d(g^{-m}(y),g^{-m}(z))=d(y,z)< c$. This is a contradiction and shows that
indeed $g$ is not elliptic.

Next assume to the contrary that $g$ is not semisimple.
By the classification of isometries of proper ${\rm CAT}(0)$-spaces 
(Proposition 3.4 of \cite{B95}), 
in this case $g$ fixes a point
$\xi\in \partial X$ and it preserves every Busemann function at $\xi$. 
Let $x\in X$ and let $H=b_\xi(x,\cdot)^{-1}(-\infty,0]$ be a closed
horoball at $\xi$. Then $H$ is a closed $g$-invariant 
convex subset of $X$ whose
closure $\overline{H}$ in $X\cup \partial X$
intersects $\partial X$ in a closed subset $\partial H$.

We claim that for $i=1,2$ 
the intersection $\partial H\cap V_i$ contains a fixed point for $g$.
For this let $\eta\in 
\partial X-\overline{V_1}-\overline{V_2}-\{\xi\}$ be
the endpoint of a contracting geodesic ray. 
Such a point exists by the
second assumption in the lemma. 
Define the shortest distance projection $\pi_H(\eta)\subset
\overline{H}$ of $\eta$ into $H$ as follows.
If an arbitrarily fixed Busemann function 
$b_\eta(x,\cdot)$ at $\eta$ assumes a 
minimum on $H$ then let $\pi_H(\eta)$ 
be the closure in $\overline{H}$ 
of the set of minima of $b_\eta(x,\cdot)$. Otherwise let 
$\pi_H(\eta)\subset \partial H$ be the 
set of accumulation points of sequences 
$(x_i)\subset H$ so that
$b_\eta(x,x_i)$ converges to the infimum of $b_\eta(x,\cdot)$ 
on $H$ 
as $i\to \infty$. Then $\pi_H(\eta)$ is a closed convex
subset of $\overline{H}$.

By Lemma \ref{visibilitypoint}, $\eta$ is a visibility point and
hence there is a geodesic $\rho$ connecting 
$\eta$ to $\xi$. This geodesic satisfies
$b_\eta(\rho(0),\rho(t))=t=-b_\xi(\rho(0),\rho(t))$ for all $t$.
In particular, if $\rho$ is parametrized in such a way
that $b_\xi(x,\rho(0))=0$ then since Busemann
functions are convex and one-Lipschitz  
we have
$\rho(0)\in \pi_H(\eta)$ (compare the simple argument in the
proof of Lemma \ref{projectionfixed}).
If $z\in \pi_{H}(\eta)\cap H$ is another point then
$z\in b_\eta(\rho(0),\cdot)^{-1}(0)\cap
b_\xi(x,\cdot)^{-1}(-\infty,0]$, and there is a 
geodesic $\rho^\prime$ connecting $\eta$ to $\xi$
which passes through $\rho^\prime(0)=z$.  
The function $t\to d(\rho(t),\rho^\prime(t))$ is
convex and bounded and hence it is constant.
By comparison, $\rho$ and $\rho^\prime$ bound a
totally geodesic embedded flat strip.

Now $\eta$ is the endpoint of a contracting
geodesic ray and hence by
Lemma \ref{thintriangle} and Lemma 3.8 of \cite{BF08},
the geodesic ray $\rho(-\infty,0]$ is $C$-contracting
for some $C>0$.
In particular, Lemma \ref{thintriangle} shows that
there is a subray of $\rho^\prime(-\infty,0]$ 
which is contained in a uniformly
bounded neighborhood of $\rho(-\infty,0]$. 
Therefore the width of the flat strip bounded 
by $\rho,\rho^\prime$ is uniformly bounded.
On the other hand, we have 
$b_\eta(\rho(0),\rho(t))=b_\eta(\rho(0),\rho^\prime(t))=-t$ 
for all $t$ and consequently the distance
between $\rho(0),\rho^\prime(0)$ is uniformly bounded as well. 
This shows that 
$\pi_{H}(\eta)\cap H$ is a bounded and hence 
compact subset of $H$.
By convexity of $\pi_{H}(\eta)$ we conclude that
$\pi_{H}(\eta)$ is contained in $H$.

The horoball $H$ is invariant under $g$ and hence
$\pi_H(g^k\eta)=g^k\pi_H(\eta)$ for all $k\in \mathbb{Z}$.
Thus if there is a compact subset of $X$ which 
intersects each of the sets $\pi_H(g^k\eta)$ $(k>0)$ 
whose diameter is uniformly bounded then $g$ has a bounded
orbit and hence $g$ is elliptic. It follows that 
there is a sequence $k_i\to \infty$ such that
$g^{k_i}\pi_H(\eta)\subset H$
leave any compact subset of $X$. 
Then up to passing to a subsequence,
the sequence $(g^{k_i}\rho(0))$ converges as $i\to \infty$ to a point
in $\partial H$. 

The horoball $H$ is closed and convex, and $g^k\rho(0)\in 
\pi_H(g^k\eta)$ for all $k$. As a consequence, for each $k\not=0$ the 
\emph{Alexandrov angle} at $g^k\rho(0)$ between
$\rho(0)$ and $g^k\eta$ is not smaller than $\pi/2$
(see \cite{BH99} for a comprehensive treatment of Alexandrov angles).
This implies that $\pi_{g^k\rho(-\infty,0]}(\rho(0))=
g^k\rho(0)$ and hence since $g^k\rho(-\infty,0]$ is 
$C$-contracting for all $k$, 
the geodesic connecting
$\rho(0)$ to $g^k\eta$ passes through a uniformly bounded neighborhood
of $g^k\rho(0)$. By comparison, this implies that the
sequence $(g^{k_i}\eta)\subset \partial X$ converges as $i\to \infty$ to 
the limit point in $\partial H$ of the sequence $(g^{k_i}\rho(0))$.
On the other hand,
for every $k$ the distance between $g^k\rho(0)$ and 
$g^{k+1}\rho(0)$ is uniformly bounded and hence
using again comparison, the limit point is a fixed
point for $g$. By the fourth assumption in the
lemma, this fixed point is contained in 
$\partial H\cap V_1$. The same argument, applied to the
sets $g^{-k}\pi_{H}(\eta)$ for $k>0$ 
shows the existence of a fixed point 
in $\partial H\cap V_2$.

By property 4) in the statement of the lemma,
this implies that the fixed point set ${\rm Fix}(g)$ of $g$  
is disconnected.
However, Corollary 3.3 of \cite{CM08} (see also
Theorem 1.1 of \cite{FNS06} for a related result) shows that
${\rm Fix}(g)$ is connected in the topology on 
$\partial X$ induced by the Tits metric $d_T$. Since the 
identity $(\partial X,d_T)\to \partial X$ is continuous
\cite{BH99}, ${\rm Fix}(g)$ is connected in the 
visual (cone) topology as well. This is a contradiction.
Together we conclude that $g$ is necessarily axial.

Let $\zeta$ be an oriented axis of $g$ and let again
$\eta\in \partial X-\overline{V_1}-\overline{V_2}-\{\xi\}$ be 
the endpoint of a contracting geodesic. 
Then $\eta$ can be connected to $\zeta(\infty)$ by a geodesic
$\rho$. This geodesic satisfies 
$b_\eta(\rho(0),\rho(t))=t$ for all $t$. Now
$t\to d(\rho(t),\zeta(t))$ is decreasing and hence by the
definition of the shortest distance projection, we conclude
that $A=\pi_{\zeta(\mathbb{R})}(\eta)$ is bounded and hence
compact.
Equivariance under the action of
$g$ implies that 
$\pi_{\zeta(\mathbb{R})}(g^k\eta)=g^kA$ for all $k\in 
\mathbb{R}$ and therefore  
$g^k\eta\to \zeta(\infty)$ $(k\to \infty)$.
On the other hand, by the fourth assumption in the lemma we have
$g^k\eta\in V_1$ for all $k\geq 1$ and hence $\zeta(\infty)\in V_1$.
The same argument also shows that $\zeta(-\infty)\in V_2$.

We are left with showing that $g$ is rank-one. Namely, we saw in 
the previous paragraph that
$g$ admits an  axis with endpoints 
$a\in V_1,b\in V_2\subset \partial X$. If $g$ is not rank-one
then $\zeta$ bounds a flat half-plane $F$ whose ideal boundary is  
an arc $\partial F$ 
connecting $a\in V_1$ to $b\in V_2$.
Then $\partial F$ intersects the open set 
$\partial X-\overline{V_1}-\overline{V_2}$.
For every $k>0$, the set $g^kF$ is a flat half-plane
bounded by $\zeta$ whose closure in $\overline{X}$ intersects
$\partial X$ in the
arc $g^k\partial F=\partial g^k F$ connecting $a$ to $b$. 

By the Arzela-Ascoli theorem, there is a sequence 
$k_i\to \infty$ such that the sequence
$(g^{k_i}F)$ of flat half-planes bounded
by $\zeta$ converges uniformly on 
compact sets to a flat half-plane $G$.
The ideal boundary $\partial G$ of $G$ intersects
$\partial X-\overline{V_1}-\overline{V_2}$ 
nontrivially. Let $z\in \partial G-\overline{V_1}-
\overline{V_2}$ be such an intersection point, let 
$x\in \zeta$ be any point and let $\beta>0$ be the
angle between $\zeta$ and the geodesic ray
$\xi$ connecting $x$ to $z=\xi(\infty)$. For sufficiently large $i$
the endpoint $z_i$ in $\partial g^{k_i}F$ of the geodesic
ray which issues from $x$, which is contained in 
$g^{k_i}F$ and which encloses an angle
$\beta$ with $\zeta$ is contained in 
$\partial X-\overline{V_1}-\overline{V_2}$.
The reasoning in the proof of Lemma \ref{northsouth}
shows that for $j>i$ the point $g^{k_j-k_i}z_i$ 
is the endpoint of the geodesic ray in $g^{k_j}F$ which
issues from $x$ and
encloses the angle $\beta$ with $\zeta$. In particular,
for $j>i$ we have $g^{k_j-k_i}z_i\in \partial X-\overline{V_1}$
which contradicts the assumption 4) in the 
statement of the lemma. 
The lemma is proven.
\end{proof}

\section{Non-elementary groups of isometries}

As in the previous sections, let $X$ be a proper ${\rm CAT}(0)$-space.
The isometry group ${\rm
Iso}(X)$ of $X$ can be equipped with a natural locally compact
$\sigma$-compact metrizable topology, the so-called \emph{compact
open topology}. With respect to this topology, a sequence
$(g_i)\subset {\rm Iso}(X)$ converges to some isometry $g$ if and
only if $g_i\to g$ uniformly on compact subsets of $X$.
A closed subset $A\subset {\rm Iso}(X)$ is compact if and only if
there is a compact subset $K$ of $X$ such that $gK\cap
K\not=\emptyset$ for every $g\in A$. In particular, the action of
${\rm Iso}(X)$ on $X$ is proper.

Let $G<{\rm Iso}(X)$ be a subgroup of the isometry group of $X$.
The \emph{limit set} $\Lambda$ of $G$ is the set
of accumulation points in $\partial X$ of one (and hence every)
orbit of the action of $G$ on $X$. If the closure of 
$G$ is non-compact then its
limit set is a compact non-empty $G$-invariant subset of $\partial X$. 
If $g\in G$ is axial with axis $\gamma$, then $\gamma(\infty),
\gamma(-\infty)\in \Lambda$. In particular, the two fixed
points for the action on $\partial X$  
of a rank-one element are contained in $\Lambda$.

A compact space is \emph{perfect} if it does not
have isolated points. 
We first observe

\begin{lemma}\label{minimal}
Let $G<{\rm Iso}(X)$ be a subgroup 
which contains a rank-one element $g$.
Then the limit set $\Lambda$ of 
$G$ is the closure in $\partial X$ of the set of 
fixed points of conjugates of $g$ in $G$.  
If $\Lambda$ contains
at least three points then 
$\Lambda$ is perfect.
\end{lemma}
\begin{proof}
Let $G<{\rm Iso}(X)$ be a subgroup which 
contains a rank-one element $g\in G$.
Let $\Lambda$ be the limit set of $G$.
We claim that $\Lambda$ is contained in the
closure of the $G$-orbit of the two fixed 
points of $g$. 
For this let $\xi\in \Lambda$, let 
$\gamma$ be a $B$-contracting axis
of $g$ for some $B>0$ and let 
$(g_i)\subset G$ 
be a sequence such that $(g_i\gamma(0))$ converges to 
$\xi$. There are two cases possible.

In the first case, up to passing to a 
subsequence, the geodesics $g_i\gamma$ eventually leave
every compact set. Let $x_0=\gamma(0)$ and for $i\geq 1$ let 
$x_i=\pi_{g_i\gamma(\mathbb{R})}(\gamma(0))$. 
Then $d(x_0,x_i)\to \infty$ $(i\to \infty)$. 
On the other hand, $g_i\gamma$ is $B$-contracting. Hence 
by Lemma \ref{thintriangle}, 
a geodesic $\zeta_i$ connecting $x_0$ to $g_ix_0$ passes through
the $3B+1$-neighborhood of $x_i$, and the same
is true for a geodesic $\eta_i$ connecting 
$x_0$ to $g_i\gamma(\infty)$.
Since $d(x_0,x_i)\to \infty$ $(i\to \infty)$, 
by convexity, 
by the description of the topology on $\partial X$ as the topology
of uniform convergence on compact sets for geodesic rays 
issuing from $x_0$, by
${\rm CAT}(0)$-comparison and compactness, we conclude the 
following.
After passing to a subsequence, the sequences
$(x_i)$ and $(g_ix_0)$ and $(g_i\gamma(\infty))$ 
converge as $i\to \infty$ to the
same point in $\partial X$. But $g_ix_0\to \xi$ and  
therefore $g_i\gamma(\infty)\to \xi$. However, $g_i\gamma(\infty)$
is a fixed point of the conjugate $g_igg_i^{-1}$ of $g$. 
This shows that
indeed $\xi$ is contained in the closure of the
fixed points of all conjugates of $g$.

In the second case there is a compact subset $K$ of $X$
such that $g_i\gamma\cap K\not=\emptyset$ for all $i$.
Since $X$ is proper by assumption, up to
passing to a subsequence we may assume that 
the $B$-contracting geodesics $g_i\gamma$
converge locally uniformly to a geodesic $\zeta$.
On the other
hand, we have $d(g_ix_0,x_0)\to \infty$ $(i\to \infty)$ 
and hence up to passing to a subsequence the
geodesic arcs connecting $x_0$ to $g_ix_0$ converge
as $i\to \infty$ to a geodesic ray which connects $x_0$ to
one of the endpoints of $\zeta$ in $\partial X$. Then 
the limit $\xi$ of the sequence $(g_ix_0)$ is an endpoint
of $\zeta$ and hence once again, $\xi$ 
is contained in the closure of the 
fixed points of conjugates of $g$ as claimed.  

Now assume that 
the limit set $\Lambda$ of $G$ contains at least 3 points. 
Let $g$ be any rank-one element of $G$. Then $\Lambda$ 
contains at least one point $\xi$ which is not a fixed point
of $g$. Since by Lemma \ref{northsouth}
$g$ acts with north-south dynamics on $\partial X$, 
the sequence $(g^k\xi)$ consists of pairwise
distinct points which converge as $k\to\infty$ 
to the attracting fixed point of $g$. Similarly,
the sequence $(g^{-k}\xi)$ converges as
$k\to\infty$ to the repelling fixed point of $g$. 
Moreover, by the above, a point $\xi\in \Lambda$ which is
not a fixed point of a rank-one element of $G$ is a limit
of fixed points of rank-one elements. This shows that
$\Lambda$ is perfect and completes the proof of the lemma.
\end{proof}

A subgroup of 
${\rm Iso}(X)$ which contains a rank-one element and 
whose limit set contains at least three points may
fix globally a point in $\partial X$. An easy example
is the upper triangular subgroup of $SL(2,\mathbb{R})$ acting
simply transitively on the hyperbolic plane ${\bf H}^2$ and fixing one
point on the boundary of ${\bf H}^2$. 
Therefore we define a 
subgroup $G$ of ${\rm Iso}(X)$ to be \emph{non-elementary} if its limit
set contains at least $3$ points and if moreover $G$ does not
fix globally a point in $\partial X$.
 
The action of a group $G$ on a topological space $Y$ is
called \emph{minimal} if every $G$-orbit is dense.  
The following lemma completes the proof of the first
part of the theorem from the introduction.

\begin{lemma}\label{moveaway}
Let $G<{\rm Iso}(X)$ be a non-elementary group 
with limit set $\Lambda$ 
which contains
a rank-one element $g\in G$ with fixed points $a\not=b\in \Lambda$.
Then for every non-empty open set $V\subset \Lambda$ 
there is some $u\in G$ with $u\{a,b\}\subset V$. 
Moreover, the action of $G$ on 
$\Lambda$ is minimal.
\end{lemma}
\begin{proof}
Let $G<{\rm Iso}(X)$ be a non-elementary subgroup
with limit set $\Lambda$ 
and let $g\in G$ be a 
rank-one isometry with attracting and 
repelling fixed points $a,b\in \Lambda$, respectively.
Let $V\subset \Lambda$ be a non-empty open set.
By Lemma \ref{minimal}, the limit set $\Lambda$ is
perfect and up to replacing 
$g$ by $g^{-1}$ (and exchanging $a$ and $b$) there is an element
$v\in G$ which maps $a$ to $v(a)\in V-\{a,b\}$. 
Then $h=vgv^{-1}$ is a
rank-one element with fixed points $v(a)\in V,v(b)
\in \Lambda$. In particular, by Lemma \ref{northsouth}, 
$h$ acts with north-south dynamics on $\Lambda$ and hence
if $v(b)\not=b$ then we have $h^k\{a,b\}\subset V$ 
for sufficiently large $k$.

If $v(b)=b$ then 
let $\rho\in G$ be an element with $\rho(b)\not=b$.
Such an element exists since by assumption, $G$ does not fix
globally the point $b$.
Since $v(a)\not=a$ 
the orbit of $a$ under the action of the infinite
cyclic subgroup $H$ of $G$ generated by $h=vgv^{-1}$ is infinite
and hence we can find
some $w\in H$ such that $w(a)\not=
\rho^{-1}(b)$. Then $u=\rho\circ w$ maps
$\{a,b\}$ to $\Lambda-\{b\}$, and for sufficiently
large $k$, the isometry $h^ku$ maps $\{a,b\}$ into $V$.

By Lemma \ref{northsouth}, a 
rank-one isometry of $X$ acts on $\partial X$ with north-south 
dynamics and hence 
every non-empty closed $G$-invariant subset $A$ of $\partial X$
contains every fixed point of every rank-one element.
Namely, if $a\not= b$ are the two fixed points of a
rank-one element of $G$ and if there is some $\xi\in A-\{a,b\}$
then also $\{a,b\}\subset A$ since $A$ is closed. On the 
other hand, if $a\in A$ then the above consideration
shows that there is some $h\in G$ with $h(a)\in \Lambda-\{a,b\}$
and once again, we conclude by invariance that $b\in A$ as well.
Now the set of all fixed points of 
rank-one elements of $G$ is $G$-invariant and hence the smallest
non-empty closed $G$-invariant subset of $\partial X$ is the
closure of the set of fixed points of rank-one elements.
This set contains
the limit set $\Lambda$ by Lemma \ref{minimal} and hence
it coincides with $\Lambda$. In other words, 
the action of $G$ on $\Lambda$ is minimal.
The lemma is proven.
\end{proof}

\begin{corollary}\label{northsouth3}
Let $G<{\rm Iso}(X)$ be a non-elementary subgroup which
contains a rank-one element. Then every element $u\in G$
which acts on $\partial X$ with north-south dynamics is rank-one.
\end{corollary}
\begin{proof}
Let $G<{\rm Iso}(X)$ be a non-elementary subgroup with limit 
set $\Lambda$ which contains a $B$-rank-one element for some $B>0$. 
Assume that $u\in G$ acts on $\partial X$ with north-south
dynamics, with fixed points $a\not=b$. 
It suffices to verify that there is some
$k>0$ and there are neighborhoods $V_1$ of $a$,
$V_2$ of $b$ which satisfy the
assumptions in Lemma \ref{northsouth2} for $u^k$. 

For this note first that since $\Lambda$ is perfect,
there is some $k>0$ and there are open neighborhoods
$V_1,V_2$ of $a,b$ in $\partial X$ with disjoint closures
$\overline{V_1},\overline{V_2}$ such that
$\Lambda-\overline{V_1}-\overline{V_2}\not=\emptyset$
and that $u^k(\partial X-V_2)\subset V_1$ and that
$u^{-k}(\partial X-V_1)\subset V_2$.
By Lemma \ref{moveaway}, there is then a biinfinite
$B$-contracting
geodesic with both endpoints in $\Lambda-\overline{V_1}-\overline{V_2}$.
This $B$-contracting geodesic is the image of an
axis of a $B$-rank-one element of $G$.

We are left with showing that via perhaps decreasing the
neighborhoods $V_1,V_2$ of $a,b$ we can guarantee that
the distance in $X$ 
between any $B$-contracting geodesic line 
with both endpoints in $V_1$
and any $B$-contracting geodesic line with both endpoints in $V_2$
is bounded 
from below by a universal positive constant. 

For this fix
a point $x\in X$ and let $\zeta_1,\zeta_2:[0,\infty)\to X$
be geodesic rays connecting $x$ to $a,b$, respectively.
By the definition of the topology on $\partial X$, 
if $U$ is a sufficiently small neighborhood of $a$ and
if $\gamma$ is a $B$-contracting geodesic line
with both endpoints
in $U$ then the geodesic connecting $x$ to 
$\pi_{\gamma(\mathbb{R})}(x)$ longes $\zeta_1$ for a 
long initial segment only depending on $U$. Similarly,
for a sufficiently small neighborhood $V$ of $b$
and any
$B$-contracting geodesic with endpoints in $V$,   
a geodesic connecting
$x$ to $\pi_{\eta(\mathbb{R})}(x)$ longes $\zeta_2$ for 
a long initial segment only depending on $V$.
From this the existence of neighborhoods $V_1,V_2$ of $a,b$ 
with the properties stated in Lemma \ref{northsouth2} is immediate.
\end{proof}

A free group with two generators is hyperbolic
in the sense of Gromov \cite{GH}. In particular,
it admits a Gromov boundary which can be viewed
as a compactification of the group. 
As another immediate consequence of 
Lemma \ref{northsouth2} we obtain
the fourth part of the theorem from the introduction.

\begin{corollary}\label{freeone}
Let $G<{\rm Iso}(X)$ be a non-elementary group which 
contains a rank-one element. Let $\Lambda\subset \partial X$
be the limit set of $G$. 
Then $G$ contains a free subgroup $\Gamma$ with two generators
and the following properties.
\begin{enumerate}
\item Every element $e\not=g\in \Gamma$ is rank-one.
\item There is a $\Gamma$-equivariant embedding of
the Gromov boundary of $\Gamma$ into $\Lambda$.
\end{enumerate}
\end{corollary}
\begin{proof}
Let $G<{\rm Iso}(X)$ be a non-elementary subgroup
which contains a rank-one element $g$.
By Lemma \ref{moveaway}, there are
two rank-one elements $g,h$ whose fixed
point sets are disjoint. 

Let $a,b$ be the attracting and
repelling fixed point of $g$, respectively,
and let $x,y$ be the attracting and
repelling fixed point of $h$.
By Lemma \ref{northsouth}, $g,h$ act with
north-south dynamics on $\partial X$.
Then up to replacing
$g,h$ by some power we can find 
small open neighborhoods 
$U_1,U_2,U_3,U_4$ of $a,
b,x,y$ in $\partial X$
with pairwise disjoint closure 
such that the pair $U_1,U_2$ satisfies the
requirements
stated in Lemma \ref{northsouth2} for $g$, and that
the same holds true for the pair $U_3,U_4$ and $h$.
Since $g,h$ act on $\partial X$ with north-south dynamics
there are numbers $k>0,\ell>0$ such that 
$g^{mk}(\overline{U_3}\cup 
\overline{U_4})\subset {U_1}\cup {U_2}$
and $h^{m\ell}(\overline{U_1}\cup \overline{U_2})\subset 
{U_3}\cup {U_4}$ 
for every $m\in \mathbb{Z}$ (here $\overline{U_i}$ is the 
closure of $U_i$). 
By the usual ping-pong argument (see e.g. p.136-138
of \cite{M88}), the isometries $g^{mk},h^{m\ell}$ generate
a free subgroup $\Gamma$ of $G$.
By construction and Lemma \ref{northsouth2}, 
the elements of $\Gamma$ 
are rank-one, and there is a 
$\Gamma$-equivariant embedding of the Gromov boundary of 
the free group into $\Lambda$.
The corollary is proven.
\end{proof}

Lemma \ref{northsouth2} is also used to show the second and third
part of the theorem from the introduction.

\begin{lemma}\label{pairfixdense}
Let $G<{\rm Iso}(X)$ be a non-elementary subgroup
which contains a rank-one element. 
\begin{enumerate}
\item 
The pairs of fixed points of rank-one elements
of $G$ are dense in $\Lambda\times \Lambda-\Delta$.
\item The action of $G$ on $\Lambda\times \Lambda-\Delta$
has a dense orbit.
\end{enumerate}
\end{lemma}
\begin{proof}
Let $G<{\rm Iso}(X)$ be a closed non-elementary
subgroup with limit set $\Lambda$ 
which contains a $B$-rank one element $g$ with
attracting fixed point $a\in \Lambda$ and repelling fixed point
$b\in \Lambda$.

Let $U\subset \Lambda\times \Lambda-\Delta$ 
be a non-empty open set. Our goal is to show
that $U$ contains a pair of fixed points
of a rank-one element $g\in G$. For this 
we may assume that there are small open sets 
$V_i\subset \partial X-\{a,b\}$ 
with disjoint closure $\overline{V_i}$ and
such that $U=V_1\times V_2\cap \Lambda\times \Lambda-\Delta$ and that 
$\Lambda-\overline{V_1}-\overline{V_2}\not=\emptyset$.
We also may assume that the distance between any $B$-contracting
geodesic with both endpoints in $\overline{V_1}$ and 
a $B$-contracting geodesic with both endpoints in $\overline{V_2}$
is bounded from below by a universal positive 
constant (compare the proof of Corollary \ref{northsouth3}).
Moreover, by Lemma \ref{moveaway}, there is some
$q\in G$ with $q\{a,b\}\subset \Lambda-\overline{V_1}-
\overline{V_2}$. The image under $q$ of a $B$-contracting axis
is a $B$-contracting geodesic with both endpoints in 
$\Lambda-\overline{V_1}-\overline{V_2}$. 

Choose some $u\in G$ which maps $\{a,b\}$ into $V_1$.
Such an element exists by Lemma \ref{moveaway}.
Then $v=ugu^{-1}$ is a $B$-rank-one isometry
with fixed points $ua,ub\in V_1$.
Similarly, there is a $B$-rank-one isometry  
$w$ with both fixed points in $V_2$. Via replacing $v,w$ by
sufficiently high powers we may assume that
$v(\partial X-V_1)\subset V_1,v^{-1}(\partial X-V_1)\subset V_1$ 
and that $w(\partial X-V_2)\subset V_2,
w^{-1}(\partial X-V_2)\subset V_2$.
Then we have $wv(\partial X-V_1)\subset V_2$ 
and $v^{-1}w^{-1}(\partial X-V_2)\subset V_1$ and hence
by Lemma \ref{northsouth2}, $wv$ is rank-one with fixed points
in $V_1\times V_2$ and hence in $U$.
The first part of the lemma is proven.

To show the second part of the lemma, we show first 
that for any non-empty open sets 
$W_1,W_2$ in $\Lambda\times \Lambda-\Delta$
there is some $h\in G$ with $hW_1\cap W_2\not=\emptyset$.
For this assume without loss of generality that
$W_1=U_1\times U_2,W_2=U_3\times U_4$ where $U_1,U_2$
and $U_3,U_4$  
are non-empty open subsets of $\Lambda$ with disjoint 
closure. Since $\Lambda$ is perfect, by possibly 
replacing $U_i$ be proper non-empty open subsets we may assume that
the sets $U_i$ are pairwise disjoint.

By the first part of the lemma, there is a rank-one element
$u\in G$ with attracting fixed point in $U_1$ and 
repelling fixed point in $U_4$.
Since $u$ acts on $\partial X$ with north-south dynamics,
there is some $k>0$ and a small open neighborhood 
$U_5\subset U_1$ of the attracting fixed point of $u$ 
such that $u^{-k}(U_5\times U_2)\subset
U_1\times U_4$. Similarly, let 
$w\in G$ be a rank-one element 
with attracting fixed point in $U_3$ and
repelling fixed point in $u^{-k}U_2\subset U_4$. 
Then we can find 
a number $\ell >0$ and an open subset 
$U_6$ of $U_2$ such that $w^\ell(u^{-k}(U_5\times U_6))\subset
U_3\times U_4$.

The limit set $\Lambda$ of $G$ 
is compact and metrizable (see \cite{BH99})
and hence $\Lambda\times \Lambda-\Delta$ is second countable.
Choose a countable basis of open sets for $\Lambda\times \Lambda-\Delta$
of the form $U_i^1\times U_i^2$ where for each $i$ the
sets $U_i^1,U_i^2$ are non-empty, open and disjoint.
We construct inductively a sequence of non-empty open sets
$V_i^j\subset U_1^j$ $(j=1,2)$ such that for each $i\geq 1$ we have 
$V_{i+1}^j\subset V_i^j$ and that there is some
$g_i\in G$ with $g_i(\overline{V_i^1}\times \overline{V_i^2})
\subset U_i^1\times U_i^2$.
Namely, write $V_1^j=U_1^j$ for $j=1,2$ and for some $i\geq 0$ 
assume that the sets $V_i^j$ and the elements
$g_i\in G$ have already been constructed. By the above,
there is some $g_{i+1}\in G$ with 
$g_{i+1}(V_i^1\times V_i^2)\cap U^1_{i+1}\times U_{i+1}^2
\not=\emptyset$. Define $V_{i+1}^1\times V_{i+1}^2$ to be any
non-empty open product set whose closure is contained in 
$V_i^1\times V_i^2\cap g^{-1}_{i+1}(U_{i+1}^1\times U_{i+1}^2)$. 
Then clearly the set $V_{i+1}^1\times V_{i+1}^2$ and the map
$g_{i+1}\in G$ satisfy the requirement of the inductive construction.

Since each of the sets $\overline{V_i^1}\times \overline{V_i^2}$ is compact
and non-empty and since 
$V_{i+1}^1\times V_{i+1}^2\subset V_i^1\times V_i^2$ for 
all $i$, there is some 
$z\in \cap_i(\overline{V_i^1}\times \overline{V_i^2})$. 
By construction, the $G$-orbit of $z$ passes through
each of the sets $U_i^1\times U_i^2$ and hence
the $G$-orbit of $z$ is dense in $\Lambda\times \Lambda-\Delta$.
This completes the proof of the lemma.
\end{proof}

We complete the discussion in this section by looking at
groups of isometries of a proper
${\rm CAT}(0)$-space $X$ which are closed with respect 
to the compact open topology.

As before, let $\Delta$ be the diagonal 
of $\partial X\times \partial X$.
We have (Lemma 6.1 of \cite{H08b}).

\begin{lemma}\label{closed}
Let $G<{\rm Iso}(X)$ be a closed subgroup with limit set $\Lambda$.
Let
$(a,b)\in \Lambda\times \Lambda-\Delta$ be the pair of fixed points of a
rank-one element. Then the $G$-orbit of $(a,b)$ is a closed
subset of $\Lambda\times \Lambda-\Delta$.
\end{lemma}

Lemma \ref{closed} is used to show 
the following strengthening of Corollary \ref{freeone}
for closed non-elementary subgroups of ${\rm Iso}(X)$ 
(compare \cite{BF02,BF08}).

\begin{proposition}\label{free}
Let $G<{\rm Iso}(X)$ be a closed non-elementary group which 
contains a rank-one element. Let $\Lambda\subset \partial X$
be the limit set of $G$. 
If $G$ does not act transitively on $\Lambda\times \Lambda-\Delta$
then $G$ contains a free subgroup $\Gamma$ with two generators
and the following properties.
\begin{enumerate}
\item Every element $e\not=g\in \Gamma$ is rank-one.
\item There is a $\Gamma$-equivariant embedding of
the Gromov boundary of $\Gamma$ into $\Lambda$.
\item  There are infinitely many
elements $u_i\in \Gamma$ $(i>0)$ with 
fixed points $a_i,b_i$
such that for all $i$ the $G$-orbit of $(a_i,b_i)
\in \Lambda\times \Lambda-\Delta$ is distinct from the orbit of
$(b_j,a_j) (j>0)$ or $(a_j,b_j)(j\not=i)$.
\end{enumerate}
\end{proposition}
\begin{proof}
Let $G<{\rm Iso}(X)$ be a closed non-elementary subgroup
with limit set $\Lambda$ which contains a rank-one element $g$
with $B$-contracting axis $\gamma$. 
Assume that $G$ does not act transitively on 
$\Lambda\times \Lambda-\Delta$. 
By Lemma \ref{pairfixdense}, there are two rank-one
elements $g,h\in G$ whose pairs of endpoints are contained
in distinct orbits for the action of $G$ on $\Lambda\times
\Lambda-\Delta$. In particular, no positive powers of these
elements are conjugate, and the elements $g,h$ admit
$B$-contracting axes for some $B>0$. By Lemma \ref{moveaway}, via replacing 
$h$ by a conjugate we may assume that the fixed
points of $g,h$ are all distinct.
Corollary \ref{freeone} then shows that up to replacing $g,h$ by
suitably chosen powers we may assume that the subgroup 
$\Gamma$ of $G$ generated by $g,h$ is free and consists
of rank-one elements. Moreover, there is an equivariant
embedding of the boundary of the free group with two generators
into the limit set of $\Lambda$ of $G$.

Now Proposition 2 of \cite{BF02} implies 
that there are infinitely
many elements in $\Gamma$ which are pairwise 
not mutually conjugate in $G$ and whose inverses
are not conjugate. We give a version of this argument here which
is in the spirit of the arguments used earlier.
Namely, let $a,b$ and $x,y$ be the attracting and repelling 
fixed points of $g,h$, respectively. 
By Lemma \ref{closed}, we may assume that
there are open neighborhoods $U_1,U_2$ of $a,b$ and
$U_3,U_4$ of $x,y$ such that the 
$G$-orbit of $(a,b)$ does not intersect
$U_3\times U_4$ and that the $G$-orbit of 
$(x,y)$ does not intersect $U_1\times U_2$.  
By replacing $g,h$ by suitable powers we may moreover
assume that 
$g(\cup_{j\not=2}
\overline{U_j})\subset U_1$,
$g^{-1}(\cup_{j\not=1}\overline{U_j})\subset U_2$,
$h(\cup_{j\not=4}\overline{U_j})\subset U_3$ and
$h^{-1}(\cup_{j\not=3}\overline{U_j})\subset U_4$.

For numbers $n,m,k,\ell>2$
consider the isometry \[f=f_{nmk\ell}=g^nh^mg^kh^{-\ell}\in \Gamma.\]
It satisfies $f(\overline{U_1})\subset U_1,
f^{-1}(\overline{U_3})\subset U_3$ and hence
the attracting fixed point of $f$ is contained in $U_1$
and its repelling fixed point is contained in $U_3$.

Since $n>2$, its conjugate
$f_1=g^{-1}fg$ satisfies $f_1(\overline{U_1})\subset U_1$
and $f_1^{-1}(\overline{U_2})\subset U_2$, i.e. its
attracting fixed point is contained in $U_1$ and its repelling
fixed point is contained in $U_2$. 
Furthermore, since $m>2$, its conjugate
$f_2=h^{-1}g^{-n}fg^nh$ has its attracting fixed point in 
$U_3$ and its repelling fixed point in $U_4$, and its conjugate
$f_3=h^{-1}fh$ has its attracting fixed point in $U_4$ and its
repelling fixed point in $U_3$.

As a consequence, $f$ is conjugate to both an element with
fixed points in $U_1\times U_2$ as well as to an
element with fixed points in $U_3\times U_4$. 
This implies that $f$ is not conjugate to either $g$ or $h$.
Moreover, since $g$ and $h$ can not both be conjugate
to $h^{-1}$, by eventually adjusting
the size of $U_3,U_4$ we may assume that
$f$ is not conjugate to $h^{-1}$.

We claim that moreover
via perhaps increasing the values of $n,\ell$ we 
can achieve that $f_{nmk\ell}$ is not conjugate to $f^{-1}$.
Namely, as $n\to \infty$, the fixed points of the conjugate
$g^{-n}f_{(2n)mkl}g^{n}=g^nh^mg^kh^{-\ell}g^n$ 
of $f_{(2n)mkl}$ converge to the fixed points of 
$g$. Similarly, the fixed points of the conjugate
$h^{-\ell}f_{nmk(2l)}^{-1}h^{\ell}=h^\ell g^{-k}h^{-m}g^{-n}h^\ell$ 
of $f_{nmk(2\ell)}^{-1}$ 
converge as $\ell\to \infty$ to the fixed points of $h$.
Thus after possibly conjugating with $g,h$, if 
$f_{nmk\ell}$ is conjugate in $G$ to $f_{nmk\ell}^{-1}$ for all $n,\ell$  
then there is a sequence of elements $g_i\in G$ which map 
a fixed compact subset $K$ of $X$ intersecting an axis for 
$g$ into a fixed compact subset $W$ of $X$ intersecting an axis for
$h$ and such that $g_i(a,b)\to (x,y)$.
However, $G$ is a closed subgroup of ${\rm Iso}(X)$ and hence
after passing to a subsequence we may assume that $g_i\to g\in G$.
Then $g(a,b)=(x,y)$ which violates the choice of 
$g,h$.

Inductively we can construct in this way a 
sequence of elements of 
$\Gamma$ with the properties stated in the proposition.
\end{proof}

\bigskip

{\bf Acknowledgement:} I am grateful to the anonymous referee
for pointing out an error in an earlier version of this paper.





\end{document}